\documentclass[11pt,a4paper]{article}
\usepackage{amsthm}
\usepackage{amsfonts}
\usepackage{amsmath,amssymb}
\usepackage{mathrsfs}
\usepackage{hyperref}
\usepackage{graphicx}
\usepackage{float}
\usepackage{indentfirst}
\newtheorem{thm}{Theorem}[section]

\newtheorem{lem}[thm]{Lemma}

\newtheorem{rem}[thm]{Remark}

\numberwithin{equation}{section}\allowdisplaybreaks
\newcommand{\s}{\,\,\,\,\,\,\,\,\,}
\begin{document}

\title{\large\bf Decay Estimations of Higher Derivatives of Solutions to Nonlinear Dirac Equation with Spin Null Structure}

\author{\normalsize \bf Zonglin Jia \\
\it Department of Mathematics and Physics, North China Electric Power University,\\
Beijing,  PR China.\\
\it Email: 50902525@ncepu.edu.cn \ \  \\
}
\date{} \maketitle
\begin{minipage}{13.5cm}
\footnotesize \bf Abstract.
We get decay rate of higher derivatives of nonlinear massless Dirac equations with a kind of ``good" spin null form. The method we rely on is similar to that of Li and Zang. However, they only give the decay rate of solution itself to nonlinear massless Dirac system.
\rm
\\

\\

\end{minipage}

\section{Introduction}\label{section1}
The conceptions and notations, that are used in the introduction, will be presented in Section \ref{section2}.

In this paper, we study Dirac equation. The famous systems are of various forms, such as linear or nonlinear massless(massive) one, ground state one, the Dirac equations coupled with other systems and so on. What we consider is defined as
\begin{eqnarray}\label{1}
\left\{ \begin{aligned}
         &\mathfrak{D}\Phi=\mathcal{N}(\Phi,\Phi)\Phi, \s\mbox{on $\mathbb{R}^{1+3}$}\\
         &\Phi(0,\cdot)=\Phi_0, \s\mbox{on $\mathbb{R}^3$}
                          \end{aligned} \right.
\end{eqnarray}
where $\Phi$ is a spin vector field, $\mathfrak{D}$ is the Dirac operator and $\mathcal{N}$ is called spin null form.
\subsection{Significance of Dirac equation}
In 1928, Paul Adrien Maurice Dirac, a British physicist, proposed a relativistic quantum mechanical equation of electron motion, namely, the Dirac equation(\cite{8}). The equation predicted the existence of antiparticles and was later confirmed by Carl Anderson's discovery of positrons in 1932.

Except for the physical meanings, this equation is also very important in math. From the point of view of partial differential equation, Dirac equation is the linearization of the famous Klein-Gordon equation. From geometric and algebraic point of view, the study of Dirac equation stimulated the development of Clifford algebra. From the point of view of complex analysis, the Dirac operator derived from the Dirac equation is a generalization of the Cauchy-Riemann operator in 4-dimensional space(\cite{1,4,9,17,18}). At present, the research on this equation mainly focuses on a single nonlinear Dirac equation and the Dirac equation coupled with other systems, such as Maxwell system.
\subsection{The previous work}
In terms of the number of equations, the Dirac equation can be divided into a single equation and coupled equations of other physical systems, for example The Dirac-Klein-Gordon equations, Maxwell-Dirac equations, etc. Its structure can be divided into linear equation, quasilinear equation and nonlinear equation. In view of the above situations, mathematicians have made in-depth research and exploration. The ground state solution of the Dirac equation(i.e. the time-independent solution) has been discussed in many literatures. For a single nonlinear equation, Bartsch and Ding proved in \cite{2} that, for any real number $p$ not less than 2, the Dirac equation has a nontrivial minimum energy solution in the $L^p$-integrable Sobolev space of order 1 when the nonlinear part satisfies some smooth condition, growth condition and periodic one. In \cite{10} Ding considers semi-classical ground state solutions of quasilinear cases, where the quasilinear part is a power function of the norm of the solution. The author proves that the equation has a minimum energy solution for all small positive parameters. And when the parameter goes to 0, the minimum energy solution of this family tends to the minimum energy solution of the related limit problem. In addition, in the process of convergence, this family of solutions will focus on the global maximum point of the potential function in a certain sense. In \cite{12}, Ding and Liu consider the case where the quasilinear part is a general function and get a conclusion similar to \cite{10}. Compared with \cite{12}, the equation considered in \cite{13} is slightly different: the quasilinear part of the equation is added with a linear function without constant terms. \cite{14} improved the conclusion of \cite{13} and its condition was weakened as: the quasilinear term is either superlinear or asymptotically linear at positive infinity; The potential function of a linear term has a local maximum but not a global maximum. In \cite{3}, Bartsch and Xu improved the result of \cite{14}, that is, the potential function of the linear term is not required to have a maximum value point, but only a critical point. For equations coupled with other physical systems, Ding and Xu studied the semi-classical ground state solutions of quasilinear Maxwell-Dirac equations and nonlinear Dirac-Klein-Gordon equations in \cite{15} and \cite{16} respectively. The quasilinear part in \cite{15} is the same as that in \cite{12}, while \cite{16} contains the nonlinear term. They extend the results of \cite{12} to equations with electromagnetic fields and equations coupled with state functions of elementary particles respectively.

Mathematicians have also made a lot of explorations on the equation of variable including time. The scattering and non-scattering problems of hartree-type nonlinear Dirac systems with critical regularity were studied by Cho, Hong and Lee in \cite{5}. The nonlinear term of this equation is not a function but a functional, and is the convolution of a function called the Yukawa potential (under some conditions it is called the Coulomb potential) with the quadratic form of the solution. In this paper, we first prove that under the Yukawa potential, if the system admits small initial data and angular regularity, the solution scatters and has global well-posedness. Through the process, the authors observe that only a small amount of angular regularity is required to obtain global well-posedness. In addition, \cite{5} studies a class of solutions with Coulomb potential are given, which are not scattering. In addition to the Dirac type equations with 3-spatial variables, the case of 2 elements has also been studied. In \cite{19}, Lee discussed two kinds of systems: one is massless honeycomb potential Dirac equation; The other is the Hartree Dirac equation. The potential of the former is a diagonal matrix of order 2, whose elements are complex quadratic forms of solutions. The nonlinear term of the latter is the convolution of the Coulomb potential with the square of the norm of the solution. In \cite{19} The results are as follows: first, these two kinds of systems have local well-posedness in square integrable fractional Sobolev space, where the order of Sobolev space is greater than a constant; second, for any square integrable Sobolev space of order less than the constant, the flow maps determined by these two kinds of systems, if them exist, are not differentiable at the origin of order 3. In the same case of 2 elements, \cite{6} considers the same Yukawa potential as \cite{5}. It is concluded that the Dirac equation is unique and exists globally with small initial data in a square integrable function space, and the solution scatterers. In order to make up for the deficiency of \cite{6}, \cite{7} studied the system of 2 elements Coulomb potential. Compared with the classical Coulomb potential, the potential function here is more general, but the conclusion is stronger: this equation has the global existence and uniqueness with small initial data in square integrable fractional Sobolev space, and the solutions scatter.

Some scholars also consider Dirac equations on manifolds. Ding and Li discussed in \cite{11} the boundary value problem of nonlinear Dirac equation on a compact Spin manifold of general dimension, which is a space-like hypersurface of a Lorentz manifold and a Riemannian manifold. The Dirac operator here is a self-adjoint elliptic operator of first order. The authors prove that the boundary value problem has a nontrivial weak solution. Further, if the nonlinear term is odd, then the problem has infinitely many weak solutions. The above result can be regarded as a discussion of the ground state solutions of Dirac equation on a general manifold. For special Lorentz manifolds, such as $(1+ N)$ dimensional Minkowski Spaces, some scholars also discuss the Dirac system above. In \cite{22}, Tzvetkov proved the uniqueness and global existence with small initial datum of a class of semilinear massless Dirac equations by contraction mapping principle, and obtained the decay properties. In \cite{20}, Li and Zang studied the asymptotic properties of linear and nonlinear Dirac equations without mass using vector field method: the solution of linear equation has negative power decay behavior with respect to spatial-time variables; Nonlinear equations (nonlinear part admits null structure) have global existence of small initial data and decay behavior of negative power, and the power is the same as linear equations.

\subsection{Our strategy}
The idea throughout this article comes from \cite{20}. Firstly, we act the Lie derivative operators of higher order with respect to vector fields set $\mathcal{A}$ to solutions of Dirac equations. Because of Theorem \ref{thm2.3} we obtain other Dirac systems with respect to higher Lie derivatives; Secondly, the Divergence Theorem in a special domain of $\mathbb{R}^{1+3}$ enables us to control a kind of energy ``$E(t)+F(t)$" which is defined in \cite{20}; Thirdly, suppose the energy is bounded by $6C_0C_1\varepsilon^2$ for $[0,T_0)$, where $C_0\varepsilon^2$ is the upper bound of $E(0)+F(0)$ and $T_0<\infty$ is the maximal time when the inequality of energy holds, we are able to change it into $4C_0C_1\varepsilon^2$ if $\varepsilon$ is sufficiently small. Thus, by the continuity of $E(t)+F(t)$ one can get a contradiction and $T_0$ must be $\infty$.

\subsection{Main result}
\begin{thm}\label{thm1}
Given $m\geq2$, under the hypothesis that $\mathcal{N}$ is $C^m$-rapidly decreasing spin null form, the nonlinear Dirac system (\ref{1}) with initial data $\Phi_0\in H^m$ admits a unique global solution
$$\Phi\in C^{m-2}(\mathbb{R}^{1+3},\mathbb{C}^4)$$
satisfying
$$
|D^k\Phi|(t,x)\lesssim\varepsilon\cdot r^{-1}\cdot\tau_{-}^{-k}\cdot\max\{r^{-1/2},\tau_{-}^{-1/2}\}\s \mbox{for}\s 0\leq k\leq m-2,
$$
provided $||\Phi_0||_{H^m}\leq\varepsilon$ with $\tau_{-}:=(1+|t-r|^2)^{1/2}$.
\end{thm}
\begin{rem}
The work what we really do is tedious computation.
\end{rem}

\section{Preliminaries and Notations}\label{section2}
In the present article, we make such appointments: Einstein summation convention are adapted; the Greek letters $\mu, \nu$ belong to $\{0,1,2,3\}$ and the Latin ones $i, j$ are in $\{1,2,3\}$, if we do not write the summation symbol ``$\sum$" .

For two quantities $Q_1$ and $Q_2$, ``$Q_1\lesssim Q_2$" means there is a universal constant $C$ such that $Q_1\leq C\cdot Q_2$, while ``$Q_1\thickapprox Q_2$" means $Q_1\lesssim Q_2$ and $Q_2\lesssim Q_1$.

Suppose $(x^0,x^1,x^2,x^3)$ is a coordinates chart of the Minkowski space $(\mathbb{R}^{1+3},\mathbf{m})$, where $\mathbf{m}$ is the Minkowski metric of signature $(-1,1,1,1)$. Sometimes we also denote $x^0$ by $t$, $(x^1,x^2,x^3)$ by $x$ and $|x|^2:=(x^1)^2+(x^2)^2+(x^3)^2$ by $r^2$.

Given a number $z\in\mathbb{C}$, $\overline{z}$ and $\mathfrak{Re}[z]$ mean the complex conjugate and real part of $z$ respectively.

\subsection{Dirac oporator}\label{subsection2.1}
The spin bundle $\mathbb{C}^4$ over Minkowski space is a trivial bundle and $D$ is its spin connection. Indeed, $D=\partial$ is the partial derivatives, since the Minkowski space is flat. Now using $D$ we can define the Dirac operator $\mathfrak{D}:=-\sqrt{-1}\gamma^{\mu}D_{\mu}$, where $\gamma^{\mu}$ are Dirac matrices which are
$$
\gamma^0:=\left(
            \begin{array}{cc}
              0 & I \\
              I & 0 \\
            \end{array}
          \right),\s
\gamma^i:=\left(
            \begin{array}{cc}
              0 & \sigma^i \\
              -\sigma^i & 0 \\
            \end{array}
          \right).
$$
Here, $I$ is the $2\times2$ unit matrix and $\sigma^i$ are Pauli's matrices:
$$
\sigma^1:=\left(
            \begin{array}{cc}
              0 & 1 \\
              1 & 0 \\
            \end{array}
          \right),\s
\sigma^2:=\left(
            \begin{array}{cc}
              0 & -\sqrt{-1} \\
              \sqrt{-1} & 0 \\
            \end{array}
          \right),\s
\sigma^3:=\left(
            \begin{array}{cc}
              1 & 0 \\
              0 & -1 \\
            \end{array}
          \right).
$$

\subsection{Spin frames}\label{subsection2.2}
Noting $\mathbb{C}^4=\mathbb{C}^2\times\mathbb{C}^2$ and for any $\Phi\in\mathbb{C}^4$, decomposing it as $\Phi:=(\psi,\chi)^T$ with $\psi,\chi\in\mathbb{C}^2$(``$T$" means the transpose), we are going to give the definitions of spin frame $\{\xi,\eta\}$ of the right-handed spin space $\mathbb{C}^2$ and the left-handed one $\{\theta,\iota\}$. Let us write $\xi:=(\xi^1,\xi^2)$, $\eta:=(\eta^1,\eta^2)$, $\theta:=(\theta_1,\theta_2)$ and $\iota:=(\iota_1,\iota_2)$, where $\xi^s,\eta^s,\theta_s, \iota_s$($s\in\{1,2\}$) are determined by the coming equations:
\begin{eqnarray*}
(|\theta_s|^2-|\iota_s|^2)\left(1-\frac{x^3}{r}\right)=\frac{x^1}{r}(\theta_s\overline{\iota}_s+\overline{\theta}_s\iota_s)+\frac{x^2}{r}(\theta_s\overline{\iota}_s-\overline{\theta}_s\iota_s)
\end{eqnarray*}
and
\begin{eqnarray*}
|\xi^s|=|\iota_s|,\s |\theta_s|=|\eta^s|,\s \theta_1\iota_2-\theta_2\iota_1=1.
\end{eqnarray*}
Hence, with respect to the new frames, $\psi$ and $\chi$ can be written as
\begin{eqnarray*}
\psi=\psi_{\xi}\cdot\xi+\psi_{\eta}\cdot\eta\s\s\mbox{and}\s\s\chi=\chi_{\theta}\cdot\theta+\chi_{\iota}\cdot\iota.
\end{eqnarray*}
The geometric interpretation of spin frames can be found in Section 2 of \cite{20}.

\subsection{Spin null form}\label{subsection2.3}
A spin null form $\mathcal{N}|_{(t,x)}:\mathbb{C}^4\times\mathbb{C}^4\rightarrow\mathbb{C}$ is a tensor field defined as follw:
\begin{eqnarray*}
\mathcal{N}|_{(t,x)}(\Phi,\Phi)&:=&A(t,x)\cdot\overline{\psi_{\xi}(t,x)}\chi_{\theta}(t,x)+B(t,x)\cdot\psi_{\xi}(t,x)\overline{\chi_{\theta}(t,x)}\\
&&+C(t,x)\cdot\overline{\psi_{\eta}(t,x)}\chi_{\iota}(t,x)+D(t,x)\cdot\psi_{\eta}(t,x)\overline{\chi_{\iota}(t,x)},
\end{eqnarray*}
where $A,B,C,D$ are complex valued functions, which are called the coefficients of the spin null form.

$\mathcal{N}$ is called $C^k$-rapidly decreasing if all of its coefficients satisfy
$$
\sup\{(t+r)^i|\partial^ju|(t,x)\s |t\geq0,\,\,\, r\geq0,\,\,\,i+j\leq k\}<\infty
$$

\subsection{Lie derivatives of spin field}\label{subsection2.4}
Given any vector field $X:=X^{\mu}\partial_{\mu}$ with $\partial_{\mu}:=\partial/\partial x^{\mu}$, the Lie derivative $\mathcal{L}_X$ is defined as
\begin{eqnarray*}
\mathcal{L}_X\Phi:=D_X\Phi-\frac{1}{4}\nabla_{[\mu,}X_{\nu]}\cdot\gamma^{\mu}\gamma^{\nu}\Phi,
\end{eqnarray*}
where $\nabla$ is the Levi-Civita connection of $(\mathbb{R}^{1+3},\mathbf{m})$ and right now $\Phi$ is a spin field or a complex valued function over the Minkowski space.

\subsection{Some special vector fields}\label{subsection2.5}
At this time we define some special vector fileds
\begin{eqnarray*}
\mathcal{A}:=\{\partial_{\mu},\,\Omega_{0i},\,\Omega_{ij},\,S\,|\,\mu=0,1,2,3,\,i,j=1,2,3\}\s\mbox{and}\s\hat{\mathcal{A}}:=\mathcal{A}\setminus\{S\},
\end{eqnarray*}
\begin{eqnarray*}
L:=(\partial_t+\partial_r)/\sqrt{2}\s\mbox{and}\s\underline{L}:=(\partial_t-\partial_r)/\sqrt{2},
\end{eqnarray*}
where
$$
\Omega_{0i}\equiv\Omega_{i0}:=x^0\partial_i+x^i\partial_0,\s\Omega_{ij}:=x^i\partial_j-x^j\partial_i,\s S:=x^{\mu}\partial_{\mu}.
$$
\begin{rem}
Denote $\omega^i:=x^i/r$, then $x^i=r\cdot\omega^i$. It is elementary to calculate $\partial_r=\omega^i\partial_i$ due to the Chain Rule.
\end{rem}

In Subsection 3.0.2 of \cite{20} we find the following relationships:
\begin{eqnarray}\label{2.0}
[\mathfrak{D},\mathcal{L}_{\partial_{\mu}}]=0,\s[\mathfrak{D},\mathcal{L}_{\Omega_{0i}}]=0,\s[\mathfrak{D},\mathcal{L}_{\Omega_{ij}}]=0,\s[\mathfrak{D},\mathcal{L}_S]=\mathfrak{D}.
\end{eqnarray}

For simplicity, let us make the coming appointment
\begin{eqnarray*}
\mathcal{L}^k_{\mathcal{A}}:=\sum\limits_{Z_1,\cdots,Z_k\in\mathcal{A}}\mathcal{L}_{Z_1}\cdots\mathcal{L}_{Z_k}\s\mbox{and}\s\mathcal{L}^k_{\hat{\mathcal{A}}}:=\sum\limits_{Z_1,\cdots,Z_k\in\hat{\mathcal{A}}}\mathcal{L}_{Z_1}\cdots\mathcal{L}_{Z_k}.
\end{eqnarray*}
For any spin field $\Phi$, $|\Phi|_{\mathcal{A},k}$ is given by
$$
|\Phi|_{\mathcal{A},k}:=\sum\limits_{i=0}^k|\mathcal{L}^i_\mathcal{A}\Phi|.
$$

Now we give the next important lemma, which can be proved by induction argument and says that the Lie derivatives acting on the spin null form keep the null structure.

\begin{lem}\label{lemma2.1}
Given a positive integer $k$, spin vector fields $\Phi,\Psi$ and a null form $\mathcal{N}$, there exist a family of null forms $\mathcal{N}_{ij}$ such that
\begin{eqnarray*}
\mathcal{L}_X^k[\mathcal{N}(\Phi,\Psi)]=\sum\limits_{i+j=0}^k\mathcal{N}_{ij}(\mathcal{L}^i_X\Phi,\mathcal{L}^j_X\Psi)
\end{eqnarray*}
for any vector field $X$. Moreover, if $X\in\mathcal{A}$ and $\mathcal{N}$ is $C^k$-rapidly decreasing, then $\mathcal{N}_{ij}$ is $C^{i+j}$-rapidly decreasing.
\end{lem}

\begin{rem}
If $\Gamma$ is a spin vector field, then it is obvious to get
\begin{eqnarray*}
\mathcal{L}_X^k[\mathcal{N}(\Phi,\Psi)\Gamma]=\sum\limits_{p+q+r=0}^k\mathcal{N}_{pq}(\mathcal{L}^p_X\Phi,\mathcal{L}^q_X\Psi)\mathcal{L}^r_X\Gamma.
\end{eqnarray*}
\end{rem}

Now we plan to compute the commutator $[\mathcal{L}^k_{\mathcal{A}},\mathfrak{\mathfrak{D}}]$. To the end, some preparation is necessary. Suppose
\begin{eqnarray}\label{2.1}
\mathcal{L}^k_{\mathcal{A}}=\mathcal{L}^{k_0}_{\hat{\mathcal{A}}}\mathcal{L}_S\mathcal{L}^{k_1}_{\hat{\mathcal{A}}}\mathcal{L}_S\cdots\mathcal{L}^{k_m}_{\hat{\mathcal{A}}}\mathcal{L}_S\mathcal{L}^{k_{m+1}}_{\hat{\mathcal{A}}}
\end{eqnarray}
with $\sum\limits_{i=0}^{m+1}k_i+m+1=k$ and $k_i\geq 0$ for all $0\leq i\leq m+1$. In the next we give the detail of our computation.

\begin{lem}\label{lemma1}
\begin{eqnarray*}
[\mathcal{L}^k_{\mathcal{A}},\mathfrak{\mathfrak{D}}]=-\sum\limits_{\substack{p+q=k-1\\p\geq0,\,q\geq0}}\mathcal{L}^p_{\mathcal{A}}\mathfrak{D}\mathcal{L}^q_{\mathcal{A}}
\end{eqnarray*}
\end{lem}
\textbf{Proof.}
From (\ref{2.0}) it follows that $[\mathcal{L}_{\hat{\mathcal{A}}},\mathfrak{D}]=0$. Hence
\begin{eqnarray*}
[\mathcal{L}^k_{\mathcal{A}},\mathfrak{\mathfrak{D}}]&=&\sum\limits_{i=0}^m\mathcal{L}^{k_0}_{\hat{\mathcal{A}}}\mathcal{L}_S\cdots\mathcal{L}^{k_i}_{\hat{\mathcal{A}}}[\mathcal{L}_S,\mathfrak{D}]\mathcal{L}^{k_{i+1}}_{\hat{\mathcal{A}}}\mathcal{L}_S\cdots\mathcal{L}^{k_m}_{\hat{\mathcal{A}}}\mathcal{L}_S\mathcal{L}^{k_{m+1}}_{\hat{\mathcal{A}}}\\
&=&\sum\limits_{\substack{p+q=k-1\\p\geq0,\,q\geq0}}\mathcal{L}^p_{\mathcal{A}}[\mathcal{L}_S,\mathfrak{D}]\mathcal{L}^q_{\mathcal{A}}=-\sum\limits_{\substack{p+q=k-1\\p\geq0,\,q\geq0}}\mathcal{L}^p_{\mathcal{A}}\mathfrak{D}\mathcal{L}^q_{\mathcal{A}}.
\end{eqnarray*}
$\hfill\qedsymbol$

Induction easily leads to the following theorem.
\begin{thm}\label{thm2.3}
There exist $c_q\in\mathbb{Z}$, the set of integers, which depends only upon $q$ and $k$, such that
\begin{eqnarray*}
\mathfrak{D}\mathcal{L}^k_{\mathcal{A}}=\mathcal{L}^k_{\mathcal{A}}\mathfrak{D}+\sum\limits_{q=0}^{k-1}c_q\cdot\mathcal{L}^q_{\mathcal{A}}\mathfrak{D}
\end{eqnarray*}
\end{thm}

Obviously, due to Theorem \ref{thm2.3} we gain that there exist a family of null forms $\mathcal{N}_{pq}$ such that
\begin{eqnarray*}
\mathfrak{D}(\mathcal{L}^k_{\mathcal{A}}\Phi)=\sum\limits_{p+q+r=0}^k\mathcal{N}_{pq}(\mathcal{L}_{\mathcal{A}}^p\Phi,\mathcal{L}^q_{\mathcal{A}}\Phi)\mathcal{L}^r_{\mathcal{A}}\Phi,
\end{eqnarray*}
if the spin vector field $\Phi$ is a solution Dirac equation.

\subsection{Some geometric domains in $\mathbb{R}^{1+3}$}\label{subsection2.6}
It is necessary to define some geometric domains, since them are used to illustrate the Divergence Theorem later.

Fix $t\in\mathbb{R}$, $\Sigma_t:=\{t\}\times\mathbb{R}^3$. If $t>0$, $\mathcal{D}_t:=\bigcup_{0\leq s\leq t}\Sigma_s$ is the space-time region between $\Sigma_0$ and $\Sigma_t$. The outgoing and incoming hypersurfaces are defined as
\begin{eqnarray*}
\mathcal{C}_u^t:=\{(s,x)|s-|x|=2u, 0\leq s\leq t\}\s \mbox{and} \s \mathcal{\hat{C}}^t_v:=\{(s,x)|s+|x|=2v,0\leq s\leq t\}.
\end{eqnarray*}
We know that the normal vector fields of $C_u^t$ and $\hat{C}^t_v$ are $\vec{n}(C_u^t)=L$ and $\vec{n}(\hat{C}_v^t)=\underline{L}$ respectively. Besides, it is also natural to define

\begin{eqnarray*}
\left\{ \begin{aligned}
         &\Sigma^+_{t,u}:=\{(t,x)|t-|x|\leq 2u\}\\
         &\Sigma^-_{t,v}:=\{(t,x)|t-|x|\geq 2v\}\\
         &\mathcal{D}^+_{t,u}:=\{(s,x)|0\leq s\leq t, s-|x|\leq 2u\}\\
         &\mathcal{D}^-_{t,v}:=\{(s,x)|0\leq s\leq t, s-|x|\geq 2v\}.
\end{aligned} \right.
\end{eqnarray*}

\subsection{Some useful norms}
We give three functions space:\\
\textbf{(a)} $H^m:=\{\phi:\mathbb{R}^3\rightarrow\mathbb{C}^4|\,\,\,\,||\phi||_{H^m}<\infty\}$
with
$$||\phi||^2_{H^m}:=\sum\limits_{i+|\alpha|=0}^m\int_{\mathbb{R}^3}(1+|x|^2)^i|\partial_x^{\alpha}\phi(x)|^2\,dx<\infty;$$
\textbf{(b)} $H^m_{\mathcal{A}}:=\{\Phi:\mathbb{R}^{1+3}\rightarrow\mathbb{C}^4|\,\,\,\,||\Phi||_{H^m_{\mathcal{A}}}<\infty\}$
with
$$||\Phi||^2_{H^m_{\mathcal{A}}}:=\sum\limits_{|k|\leq m,Z_1,\cdots,Z_k\in\mathcal{A}}||Z_1\cdots Z_k\Phi||^2_{L^2};$$
\textbf{(c)} $H^m_{\mathcal{L}_{\mathcal{A}}}:=\{\Phi:\mathbb{R}^{1+3}\rightarrow\mathbb{C}^4|\,\,\,\,||\Phi||_{H^m_{\mathcal{L}_{\mathcal{A}}}}<\infty\}$
with
$$||\Phi||^2_{H^m_{\mathcal{L}_{\mathcal{A}}}}:=\sum\limits_{|k|\leq m,Z_1,\cdots,Z_k\in\mathcal{A}}||\mathcal{L}_{Z_1}\cdots \mathcal{L}_{Z_k}\Phi||^2_{L^2};$$
In the proof of Lemma 4.0.1 of \cite{21} we find that ``$||\cdot||_{H^m_{\mathcal{A}}}$" is equivalent to ``$||\cdot||_{H^m_{\mathcal{L}_{\mathcal{A}}}}$".
\section{Proof of The Main Theorem}\label{section3}
We split the Dirac equation (\ref{1}) into two equations
\begin{eqnarray*}
\delta^{\mu}D_{\mu}\chi=\sqrt{-1}\mathcal{N}(\Phi,\Phi)\psi\s\mbox{and}\s\hat{\delta}^{\mu}D_{\mu}\psi=\sqrt{-1}\mathcal{N}(\Phi,\Phi)\chi
\end{eqnarray*}
with
\begin{eqnarray*}
\delta^0=\hat{\delta}^0=I\s\mbox{and}\s \delta^j=-\hat{\delta}^j=\sigma^j.
\end{eqnarray*}
\subsection{Auxiliary tensor}\label{subsection3.1}
Define two auxiliary tensors
\begin{eqnarray*}
A^{\mu}(\psi):=\psi^{\dag}\hat{\delta}^{\mu}\psi \s \mbox{and} \s B^{\mu}(\chi):=\chi^{\dag}\delta^{\mu}\chi,
\end{eqnarray*}
where ``$\dag$" is the complex conjugate transpose. It is easy to see
\begin{eqnarray*}
\nabla_{\mu}A^{\mu}(\psi)=2\mathfrak{Re}\left[\sqrt{-1}\mathcal{N}(\Phi,\Phi)\psi^{\dag}\chi\right] \s \mbox{and} \s  \nabla_{\mu}B^{\mu}(\chi)=2\mathfrak{Re}\left[\sqrt{-1}\mathcal{N}(\Phi,\Phi)\chi^{\dag}\psi\right].
\end{eqnarray*}
Employing the same method we also obtain the next two identities:
\begin{eqnarray}\label{4}
\nabla_{\mu}A^{\mu}(\mathcal{L}^k_{\mathcal{A}}\psi)=\sum\limits_{p+q+r=0}^k\mathfrak{Re}\left[\sqrt{-1}\mathcal{N}_{pq}(\mathcal{L}^p_{\mathcal{A}}\Phi,\mathcal{L}^q_{\mathcal{A}}\Phi)(\mathcal{L}^k_{\mathcal{A}}\psi)^{\dag}\mathcal{L}^r_{\mathcal{A}}\chi\right]
\end{eqnarray}
and
\begin{eqnarray*}
\nabla_{\mu}B^{\mu}(\mathcal{L}^k_{\mathcal{A}}\chi)=\sum\limits_{p+q+r=0}^k\mathfrak{Re}\left[\sqrt{-1}\mathcal{N}_{pq}(\mathcal{L}^p_{\mathcal{A}}\Phi,\mathcal{L}^q_{\mathcal{A}}\Phi)(\mathcal{L}^k_{\mathcal{A}}\chi)^{\dag}\mathcal{L}^r_{\mathcal{A}}\psi\right].
\end{eqnarray*}

\subsection{Uniqueness and global existence with small initial data}\label{subsection3.2}
Referred to \cite{22}, in a suitable space Tzvetkov gives the uniqueness and global existence with small initial data to nonlinear massless Dirac system of the following form
\begin{eqnarray}\label{2}
\left\{ \begin{aligned}
         &\mathfrak{D}\Phi=F(\Phi), \s\mbox{on $\mathbb{R}^{1+3}$}\\
         &\Phi(0,\cdot)=\Phi_0, \s\mbox{on $\mathbb{R}^3$}
                          \end{aligned} \right.
\end{eqnarray}
with $|F(\Phi)|\lesssim|\Phi|^p$ for $p\geq2$. His idea is contraction mapping principle and the chosen space is $X^m_{\epsilon}:=\{\Phi:\mathbb{R}^{1+3}\rightarrow\mathbb{C}^4|\,\,\,\,||\Phi||_{H^m_{\mathcal{A}}}\leq\epsilon\}(m=2)$. Moreover, the set which the initial data belongs to is $Y^m_{\eta}:=\{\phi:\mathbb{R}^3\rightarrow\mathbb{C}^4|\,\,\,\,||\phi||_{H^m}\leq\eta\}(m=2)$. As long as modify his proof a little, we can gain the following existence result.

\begin{thm}\label{thm3.1}
Given $m\geq2$, there exists an $a_0>0$ such that for all $a\in(0,a_0]$, if $\Phi_0\in Y^m_{a}$, (\ref{1}) admits a unique global solution $\Phi\in X^m_{\epsilon}$, provided $\epsilon$ is sufficiently small.
\end{thm}

\begin{rem}
Like \cite{20}, \cite{22} gives the decay estimations of solutions to (\ref{2}). Therefore, what we will do is to estimate the decay rate of higher derivatives.
\end{rem}

\subsection{Decay estimations of higher derivatives}\label{subsection3.3}
The Divergence Theorem yields
\begin{eqnarray}\label{3}
&&\int_{\Sigma^+_{t,u}}A^{\mu}(\mathcal{L}^k_{\mathcal{A}}\psi)\cdot\vec{n}_{\mu}(\Sigma^+_{t,u})+\int_{\Sigma^+_{0,u}}A^{\mu}(\mathcal{L}^k_{\mathcal{A}}\psi)\cdot\vec{n}_{\mu}(\Sigma^+_{0,u})\nonumber\\
&+&\int_{\mathcal{C}^t_u}A^{\mu}(\mathcal{L}^k_{\mathcal{A}}\psi)\cdot\vec{n}_{\mu}(\mathcal{C}^t_u)=\iint\limits_{\mathcal{D}^+_{t,u}}\nabla_{\mu}A^{\mu}(\mathcal{L}^k_{\mathcal{A}}\psi).
\end{eqnarray}
Substituting (\ref{4}) into (\ref{3}) leads to
\begin{eqnarray}\label{5}
&&\int_{\Sigma^+_{t,u}}A^{\mu}(\mathcal{L}^k_{\mathcal{A}}\psi)\cdot\vec{n}_{\mu}(\Sigma^+_{t,u})+\int_{\Sigma^+_{0,u}}A^{\mu}(\mathcal{L}^k_{\mathcal{A}}\psi)\cdot\vec{n}_{\mu}(\Sigma^+_{0,u})+\int_{\mathcal{C}^t_u}A^{\mu}(\mathcal{L}^k_{\mathcal{A}}\psi)\cdot\vec{n}_{\mu}(\mathcal{C}^t_u)\nonumber\\
&=&\sum\limits_{p+q+r=0}^k\iint\limits_{\mathcal{D}^+_{t,u}}\mathfrak{Re}\left[\sqrt{-1}\mathcal{N}_{pq}(\mathcal{L}^p_{\mathcal{A}}\Phi,\mathcal{L}^q_{\mathcal{A}}\Phi)(\mathcal{L}^k_{\mathcal{A}}\psi)^{\dag}\mathcal{L}^r_{\mathcal{A}}\chi\right].
\end{eqnarray}
According to the Infeld-van der Waerden symbols which is stated amply in Lemma 1.6.1 of \cite{20}, we gain
\begin{eqnarray*}
\int_{\Sigma^+_{t,u}}|\mathcal{L}^k_{\mathcal{A}}\psi|^2+\int_{\mathcal{C}^t_u}|(\mathcal{L}^k_{\mathcal{A}}\psi)_{\eta}|^2&\lesssim&\int_{\Sigma^+_{0,u}}|\mathcal{L}^k_{\mathcal{A}}\psi|^2+\sum\limits_{p+q+r=0}^k\iint\limits_{\mathcal{D}^+_{t,u}}|\mathcal{L}^p_{\mathcal{A}}\psi|\cdot|\mathcal{L}^q_{\mathcal{A}}\chi|\cdot|\mathcal{L}^r_{\mathcal{A}}\chi|\cdot|\mathcal{L}^k_{\mathcal{A}}\psi|\\
&\leq&\int_{\Sigma_0}|\mathcal{L}^k_{\mathcal{A}}\psi|^2+\sum\limits_{p+q+r=0}^k\iint\limits_{\mathcal{D}_t}|\mathcal{L}^p_{\mathcal{A}}\psi|\cdot|\mathcal{L}^q_{\mathcal{A}}\chi|\cdot|\mathcal{L}^r_{\mathcal{A}}\chi|\cdot|\mathcal{L}^k_{\mathcal{A}}\psi|
\end{eqnarray*}
where we have used the fact that $\mathcal{N}_{pq}$ is $C^{p+q}$-rapidly decreasing which can be deduced from the hypothesis in Theorem \ref{thm1}. Similarly, the following energy estimation holds true
\begin{eqnarray*}
\int_{\Sigma^-_{t,v}}|\mathcal{L}^k_{\mathcal{A}}\psi|^2+\int_{\hat{\mathcal{C}}^t_v}|(\mathcal{L}^k_{\mathcal{A}}\psi)_{\xi}|^2&\lesssim&\int_{\Sigma^-_{0,v}}|\mathcal{L}^k_{\mathcal{A}}\psi|^2+\sum\limits_{p+q+r=0}^k\iint\limits_{\mathcal{D}^-_{t,v}}|\mathcal{L}^p_{\mathcal{A}}\psi|\cdot|\mathcal{L}^q_{\mathcal{A}}\chi|\cdot|\mathcal{L}^r_{\mathcal{A}}\chi|\cdot|\mathcal{L}^k_{\mathcal{A}}\psi|\\
&\leq&\int_{\Sigma_0}|\mathcal{L}^k_{\mathcal{A}}\psi|^2+\sum\limits_{p+q+r=0}^k\iint\limits_{\mathcal{D}_t}|\mathcal{L}^p_{\mathcal{A}}\psi|\cdot|\mathcal{L}^q_{\mathcal{A}}\chi|\cdot|\mathcal{L}^r_{\mathcal{A}}\chi|\cdot|\mathcal{L}^k_{\mathcal{A}}\psi|.
\end{eqnarray*}
Adding the above two inequalities gives
\begin{eqnarray*}
\int_{\Sigma_t}|\mathcal{L}^k_{\mathcal{A}}\psi|^2+\sup\limits_{u\in\mathbb{R}}\int_{\mathcal{C}^t_u}|(\mathcal{L}^k_{\mathcal{A}}\psi)_{\eta}|^2+\sup\limits_{v\in\mathbb{R}}\int_{\hat{\mathcal{C}}^t_v}|(\mathcal{L}^k_{\mathcal{A}}\psi)_{\xi}|^2\lesssim\int_{\Sigma_0}|\mathcal{L}^k_{\mathcal{A}}\psi|^2+\sum\limits_{p+q+r=0}^k\mathbf{I(p,q,r,k)}
\end{eqnarray*}
with
\begin{eqnarray*}
\mathbf{I(p,q,r,k)}:=\iint\limits_{\mathcal{D}_t}|\mathcal{L}^p_{\mathcal{A}}\psi|\cdot|\mathcal{L}^q_{\mathcal{A}}\chi|\cdot|\mathcal{L}^r_{\mathcal{A}}\chi|\cdot|\mathcal{L}^k_{\mathcal{A}}\psi|.
\end{eqnarray*}
The same method can be applied to $\mathcal{L}^k_{\mathcal{A}}\chi$ and $B^{\mu}(\mathcal{L}^k_{\mathcal{A}}\chi)$ in order to get
\begin{eqnarray*}
\int_{\Sigma_t}|\mathcal{L}^k_{\mathcal{A}}\chi|^2+\sup\limits_{u\in\mathbb{R}}\int_{\mathcal{C}^t_u}|(\mathcal{L}^k_{\mathcal{A}}\chi)_{\theta}|^2+\sup\limits_{v\in\mathbb{R}}\int_{\hat{\mathcal{C}}^t_v}|(\mathcal{L}^k_{\mathcal{A}}\chi)_{\iota}|^2\lesssim\int_{\Sigma_0}|\mathcal{L}^k_{\mathcal{A}}\chi|^2+\sum\limits_{p+q+r=0}^k\mathbf{II(p,q,r,k)}
\end{eqnarray*}
with
\begin{eqnarray*}
\mathbf{II(p,q,r,k)}:=\iint\limits_{\mathcal{D}_t}|\mathcal{L}^p_{\mathcal{A}}\psi|\cdot|\mathcal{L}^q_{\mathcal{A}}\chi|\cdot|\mathcal{L}^r_{\mathcal{A}}\psi|\cdot|\mathcal{L}^k_{\mathcal{A}}\chi|.
\end{eqnarray*}

Define
\begin{eqnarray*}
E(t):=\sum\limits_{k=0}^m\int_{\Sigma_t}\left[|\mathcal{L}^k_{\mathcal{A}}\psi|^2+|\mathcal{L}^k_{\mathcal{A}}\chi|^2\right]
\end{eqnarray*}
and
\begin{eqnarray*}
F(t):=\sum\limits_{k=0}^m\sup\limits_{u\in\mathbb{R}}\int_{\mathcal{C}^t_u}\left[|(\mathcal{L}^k_{\mathcal{A}}\psi)_{\eta}|^2+|(\mathcal{L}^k_{\mathcal{A}}\chi)_{\theta}|^2\right]+\sum\limits_{k=0}^m\sup\limits_{v\in\mathbb{R}}\int_{\hat{\mathcal{C}}^t_v}\left[|(\mathcal{L}^k_{\mathcal{A}}\psi)_{\xi}|^2+|(\mathcal{L}^k_{\mathcal{A}}\chi)_{\iota}|^2\right].
\end{eqnarray*}
At this time let us employ continuity argument. Suppose for all $t\in[0,T_0)$, the following estimation holds
\begin{eqnarray}\label{6}
E(t)+F(t)\leq 6C_0\cdot C_1\cdot\varepsilon^2,
\end{eqnarray}
where $C_0$ is a constant with
$$
E(0)\leq C_0\cdot\varepsilon^2
$$
and $T_0$ is the maximal time such that (\ref{6}) holds. Our goal is to prove that there is a positive number $\varepsilon_0$, which is universal, such that for all $\varepsilon\in(0,\varepsilon_0]$, the coming inequality is true
$$
E(t)+F(t)\leq 4C_0\cdot C_1\cdot\varepsilon^2\s\s\mbox{for all $t\in[0,T_0)$}.
$$
Therefore, $T_0$ can be extended to $\infty$. The key of our proof is estimating
$$
\sum\limits_{p+q+r=0}^m\left\{\mathbf{I(p,q,r,m)}+\mathbf{II(p,q,r,m)}\right\}.
$$
We only deal with $\sum\limits_{p+q+r=0}^m\mathbf{I(p,q,r,m)}$ since the other terms can be treated by the same way.

Lemma 4.3.2 and Lemma 4.3.5 of \cite{20} yield
\begin{eqnarray}\label{7}
|\phi|^2\lesssim r^{-3}\cdot\sum\limits_{i=0}^2\int_{\mathcal{C}_u^t}|\mathcal{L}^i_{\mathcal{A}}\phi|^2\s\mbox{for $t>0$}
\end{eqnarray}
and
\begin{eqnarray}\label{8}
|\phi|^2\lesssim r^{-2}\tau_{-}^{-1}\cdot\sum\limits_{i=0}^2\int_{\hat{\mathcal{C}}_v^t}|\mathcal{L}^i_{\mathcal{A}}\phi|^2\s\mbox{for $v>0$.}
\end{eqnarray}

Given $X\in\mathcal{A}$, from Lemma 3.0.2 of \cite{20} it follows that
\begin{eqnarray*}
\mathcal{L}_X\xi=\frac{1}{4}(\nabla_{\mu}X_{\nu}-\nabla_{\nu}X_{\mu})\hat{\delta}^{\mu}\hat{\delta}^{\nu}\cdot\xi.
\end{eqnarray*}
Because the coefficients of $X$ are all linear functions, $(\nabla_{\mu}X_{\nu}-\nabla_{\nu}X_{\mu})\hat{\delta}^{\mu}\hat{\delta}^{\nu}/4$ is a constant. Induction and Proposition 3.0.3 of \cite{20} lead to the following lemma.

\begin{lem}\label{lem3.1}
There are constants sets $\{a_i,\tilde{a}_i|0\leq i\leq k-1\}$, $\{b_i,\tilde{b}_i|0\leq i\leq k-1\}$, $\{c_i,\tilde{c}_i|0\leq i\leq k-1\}$ and $\{d_i,\tilde{d}_i|0\leq i\leq k-1\}$ satisfying
\begin{eqnarray*}
&&(\mathcal{L}^k_X\psi)_{\xi}=\mathcal{L}^k_X\psi_{\xi}+\sum\limits_{i=0}^{k-1}a_i\cdot\mathcal{L}^i_X\psi_{\xi},\s\mathcal{L}^k_X\psi_{\xi}=(\mathcal{L}^k_X\psi)_{\xi}+\sum\limits_{i=0}^{k-1}\tilde{a}_i\cdot(\mathcal{L}^i_X\psi)_{\xi},\\
&&(\mathcal{L}^k_X\psi)_{\eta}=\mathcal{L}^k_X\psi_{\eta}+\sum\limits_{i=0}^{k-1}b_i\cdot\mathcal{L}^i_X\psi_{\eta},\s\mathcal{L}^k_X\psi_{\eta}=(\mathcal{L}^k_X\psi)_{\eta}+\sum\limits_{i=0}^{k-1}\tilde{b}_i\cdot(\mathcal{L}^i_X\psi)_{\eta},\\
&&(\mathcal{L}^k_X\chi)_{\theta}=\mathcal{L}^k_X\chi_{\theta}+\sum\limits_{i=0}^{k-1}c_i\cdot\mathcal{L}^i_X\chi_{\theta},\s\mathcal{L}^k_X\chi_{\theta}=(\mathcal{L}^k_X\chi)_{\theta}+\sum\limits_{i=0}^{k-1}\tilde{c}_i\cdot(\mathcal{L}^i_X\chi)_{\theta},\\
&&(\mathcal{L}^k_X\chi)_{\iota}=\mathcal{L}^k_X\chi_{\iota}+\sum\limits_{i=0}^{k-1}d_i\cdot\mathcal{L}^i_X\chi_{\iota},\s\mathcal{L}^k_X\chi_{\iota}=(\mathcal{L}^k_X\chi)_{\iota}+\sum\limits_{i=0}^{k-1}\tilde{d}_i\cdot(\mathcal{L}^i_X\chi)_{\iota}
\end{eqnarray*}
provided $X\in\mathcal{A}$. Therefore, the readers can obtain
\begin{eqnarray*}
|\psi|_{\mathcal{A},k}\thickapprox|\psi_{\xi}|_{\mathcal{A},k}+|\psi_{\eta}|_{\mathcal{A},k}\s\mbox{and}\s|\chi|_{\mathcal{A},k}\thickapprox|\chi_{\theta}|_{\mathcal{A},k}+|\chi_{\iota}|_{\mathcal{A},k}.
\end{eqnarray*}
\end{lem}
Thanks to Lemma \ref{lem3.1}, we claim
$$
F(t)\thickapprox\sum\limits_{k=0}^m\sup\limits_{u\in\mathbb{R}}\int_{\mathcal{C}^t_u}\left[|\mathcal{L}^k_{\mathcal{A}}\psi_{\eta}|^2+|\mathcal{L}^k_{\mathcal{A}}\chi_{\theta}|^2\right]+\sum\limits_{k=0}^m\sup\limits_{v\in\mathbb{R}}\int_{\hat{\mathcal{C}}^t_v}\left[|\mathcal{L}^k_{\mathcal{A}}\psi_{\xi}|^2+|\mathcal{L}^k_{\mathcal{A}}\chi_{\iota}|^2\right].
$$
It is easy to see
$$F(t)\leq6C_0\cdot C_1\cdot\varepsilon^2$$
implies
$$
\sum\limits_{k=0}^m\int_{\mathcal{C}^t_{-1/2}}\left[|\mathcal{L}^k_{\mathcal{A}}\psi_{\eta}|^2+|\mathcal{L}^k_{\mathcal{A}}\chi_{\theta}|^2\right]+\sum\limits_{k=0}^m\int_{\hat{\mathcal{C}}^t_{t+1/2}}\left[|\mathcal{L}^k_{\mathcal{A}}\psi_{\xi}|^2+|\mathcal{L}^k_{\mathcal{A}}\chi_{\iota}|^2\right]\lesssim\varepsilon^2.
$$
Note that on $\mathcal{C}^t_{-1/2}$, one can easily get $r=t+1$; similarly, on $\hat{\mathcal{C}}^t_{t+1/2}$, there holds $r\geq t+1$ due to the definition of $\hat{\mathcal{C}}^t_v$. Recall that (\ref{7}) and (\ref{8}) imply
$$
\sum\limits_{k=0}^{m-2}\left[|\mathcal{L}^k_{\mathcal{A}}\psi_{\eta}|^2+|\mathcal{L}^k_{\mathcal{A}}\chi_{\theta}|^2\right]\lesssim\varepsilon^2\cdot r^{-3}
$$
and
$$
\sum\limits_{k=0}^{m-2}\left[|\mathcal{L}^k_{\mathcal{A}}\psi_{\xi}|^2+|\mathcal{L}^k_{\mathcal{A}}\chi_{\iota}|^2\right]\lesssim\varepsilon^2\cdot\tau_{-}^{-1}\cdot r^{-2}
$$
respectively. Hence
$$
\sum\limits_{k=0}^{m-2}\left[|\mathcal{L}^k_{\mathcal{A}}\psi_{\eta}|^2+|\mathcal{L}^k_{\mathcal{A}}\chi_{\theta}|^2\right]\lesssim\varepsilon^2\cdot (t+1)^{-3}
$$
and
$$
\sum\limits_{k=0}^{m-2}\left[|\mathcal{L}^k_{\mathcal{A}}\psi_{\xi}|^2+|\mathcal{L}^k_{\mathcal{A}}\chi_{\iota}|^2\right]\lesssim\varepsilon^2\cdot (t+1)^{-2},
$$
which are equivalent to
$$
\sum\limits_{k=0}^{m-2}\left[|(\mathcal{L}^k_{\mathcal{A}}\psi)_{\eta}|^2+|(\mathcal{L}^k_{\mathcal{A}}\chi)_{\theta}|^2\right]\lesssim\varepsilon^2\cdot (t+1)^{-3}
$$
and
$$
\sum\limits_{k=0}^{m-2}\left[|(\mathcal{L}^k_{\mathcal{A}}\psi)_{\xi}|^2+|(\mathcal{L}^k_{\mathcal{A}}\chi)_{\iota}|^2\right]\lesssim\varepsilon^2\cdot (t+1)^{-2}.
$$

Now we estimate $\sum\limits_{p+q+r=0}^m\mathbf{I(p,q,r,m)}$. This quantity can be bounded by
\begin{eqnarray*}
&&\sum\limits_{p+q+r=0}^m\int_0^t\,ds\left\{\int_{\Sigma_s}|\mathcal{L}^p_{\mathcal{A}}\psi|^2\cdot|\mathcal{L}^q_{\mathcal{A}}\chi|^2\cdot|\mathcal{L}^r_{\mathcal{A}}\chi|^2\right\}^{1/2}\left\{\int_{\Sigma_s}|\mathcal{L}^m_{\mathcal{A}}\psi|^2\right\}^{1/2}\\
&\lesssim&\sum\limits_{p+q+r=0}^m\int_0^t\,ds\left\{\int_{\Sigma_s}|\mathcal{L}^p_{\mathcal{A}}\psi|^2\cdot|\mathcal{L}^q_{\mathcal{A}}\chi|^2\cdot|\mathcal{L}^r_{\mathcal{A}}\chi|^2\right\}^{1/2}E(s)^{1/2}\\
&\lesssim&\varepsilon\sum\limits_{p+q+r=0}^m\int_0^t\,ds\left\{\int_{\Sigma_s}|\mathcal{L}^p_{\mathcal{A}}\psi|^2\cdot|\mathcal{L}^q_{\mathcal{A}}\chi|^2\cdot|\mathcal{L}^r_{\mathcal{A}}\chi|^2\right\}^{1/2}\\
&=&\varepsilon\left(\sum\limits_{\substack{p+q+r=0\\p\leq m-2}}^m+\sum\limits_{\substack{p+q+r=0\\p= m-1}}^m+\sum\limits_{\substack{p+q+r=0\\p= m}}^m\right)\int_0^t\,ds\left\{\int_{\Sigma_s}|\mathcal{L}^p_{\mathcal{A}}\psi|^2\cdot|\mathcal{L}^q_{\mathcal{A}}\chi|^2\cdot|\mathcal{L}^r_{\mathcal{A}}\chi|^2\right\}^{1/2}\\
&\lesssim&\varepsilon^2\sum\limits_{\substack{p+q+r=0\\p\leq m-2}}^m\int_0^t\,ds\left\{\int_{\Sigma_s}|\mathcal{L}^q_{\mathcal{A}}\chi|^2\cdot|\mathcal{L}^r_{\mathcal{A}}\chi|^2\right\}^{1/2}(s+1)^{-1}+\varepsilon\int_0^t\,ds\left\{\int_{\Sigma_s}|\chi|^4\cdot|\mathcal{L}^m_{\mathcal{A}}\psi|^2\right\}^{1/2}\\
&&+\varepsilon\sum\limits_{q+r=0}^1\int_0^t\,ds\left\{\int_{\Sigma_s}|\mathcal{L}^q_{\mathcal{A}}\chi|^2\cdot|\mathcal{L}^r_{\mathcal{A}}\chi|^2\cdot|\mathcal{L}^{m-1}_{\mathcal{A}}\psi|^2\right\}^{1/2}:=\varepsilon^2\cdot\mathbf{I}_1+\varepsilon\cdot\mathbf{I}_2+\varepsilon\cdot\mathbf{I}_3.
\end{eqnarray*}

Note that
\begin{eqnarray*}
\mathbf{I}_3&\lesssim&\varepsilon^2\int_0^t\,ds\left\{\int_{\Sigma_s}|\mathcal{L}^{m-1}_{\mathcal{A}}\psi|^2\right\}^{1/2}\cdot(s+1)^{-2}\leq\varepsilon^2\int_0^t(s+1)^{-2}E(s)^{1/2}\,ds\\
&\lesssim&\varepsilon^3\int_0^t(s+1)^{-2}\,ds\leq\varepsilon^3
\end{eqnarray*}
and
\begin{eqnarray*}
\mathbf{I}_2\lesssim\varepsilon^2\int_0^t(s+1)^{-2}E(s)^{1/2}\,ds\lesssim\varepsilon^3\int_0^t(s+1)^{-2}\,ds\leq\varepsilon^3.
\end{eqnarray*}
By some obvious observation we have
\begin{eqnarray*}
\mathbf{I}_1&\leq&2\sum\limits_{\substack{p+q+r=0\\p\leq m-2,\,\,q\leq r}}^m\int_0^t\,ds\left\{\int_{\Sigma_s}|\mathcal{L}^q_{\mathcal{A}}\chi|^2\cdot|\mathcal{L}^r_{\mathcal{A}}\chi|^2\right\}^{1/2}(s+1)^{-1}\\
&\leq&2\sum\limits_{\substack{p+q+r=0\\p\leq m-2,\,\,q\leq m/2}}^m\int_0^t\,ds\left\{\int_{\Sigma_s}|\mathcal{L}^q_{\mathcal{A}}\chi|^2\cdot|\mathcal{L}^r_{\mathcal{A}}\chi|^2\right\}^{1/2}(s+1)^{-1}\\
&\lesssim&\varepsilon\sum\limits_{r=0}^m\int_0^t\,ds\left\{\int_{\Sigma_s}|\mathcal{L}^r_{\mathcal{A}}\chi|^2\right\}^{1/2}(s+1)^{-2}\lesssim\varepsilon\int_0^t\,ds\,\,E(s)^{1/2}\cdot(s+1)^{-2}\\
&\lesssim&\varepsilon^2\int_0^t\,ds\,\,\cdot(s+1)^{-2}\leq\varepsilon^2.
\end{eqnarray*}
In conclusion, we obtain
$$
\sum\limits_{p+q+r=0}^m\mathbf{I(p,q,r,m)}\lesssim\varepsilon^4
$$
implying
$$
E(t)+F(t)\leq C_0\cdot\varepsilon^2+C_6\cdot\varepsilon^4\s\s\mbox{for all $t\in[0,T_0)$}.
$$
Pick $\varepsilon_0:=\{(4C_1-1)\cdot C_0/C_6\}^{1/2}$. Given $\varepsilon\in(0,\varepsilon_0]$, it is not difficult to have
$$
E(t)+F(t)\leq 4C_0\cdot C_1\cdot\varepsilon^2\s\s\mbox{for all $t\in[0,T_0)$}.
$$
Hence we get a contradiction and $T_0$ must be $\infty$.

From Lemma 4.3.1 of \cite{20} it follows that
$$
|D^k\Phi|\lesssim\tau_{-}^{-k}\cdot|\Phi|_{\mathcal{A},k}.
$$
Lemma 4.3.5 of \cite{20} yieds
\begin{eqnarray*}
|\mathcal{L}^k_{\mathcal{A}}\psi_{\eta}|\lesssim r^{-3/2}\cdot\left\{\sum\limits_{i=0}^2\int_{\mathcal{C}^t_u}|\mathcal{L}^{k+i}_{\mathcal{A}}\psi_{\eta}|^2\right\}^{1/2}\mbox{and} \s|\mathcal{L}^k_{\mathcal{A}}\chi_{\theta}|\lesssim r^{-3/2}\cdot\left\{\sum\limits_{i=0}^2\int_{\mathcal{C}^t_u}|\mathcal{L}^{k+i}_{\mathcal{A}}\chi_{\theta}|^2\right\}^{1/2},
\end{eqnarray*}
\begin{eqnarray*}
|\mathcal{L}^k_{\mathcal{A}}\psi_{\xi}|\lesssim r^{-1}\tau_{-}^{-1/2}\cdot\left\{\int_{\hat{\mathcal{C}}^t_v}\left(\sum\limits_{i=0}^2|\mathcal{L}^{k+i}_{\mathcal{A}}\psi_{\xi}|^2+\tau_{-}^2\sum\limits_{j=0}^1|\mathcal{L}_{\underline{L}}\mathcal{L}^{k+j}_{\mathcal{A}}\psi_{\xi}|^2\right)\right\}^{1/2}
\end{eqnarray*}
and
\begin{eqnarray*}
|\mathcal{L}^k_{\mathcal{A}}\chi_{\iota}|\lesssim r^{-1}\tau_{-}^{-1/2}\cdot\left\{\int_{\hat{\mathcal{C}}^t_v}\left(\sum\limits_{i=0}^2|\mathcal{L}^{k+i}_{\mathcal{A}}\chi_{\iota}|^2+\tau_{-}^2\sum\limits_{j=0}^1|\mathcal{L}_{\underline{L}}\mathcal{L}^{k+j}_{\mathcal{A}}\chi_{\iota}|^2\right)\right\}^{1/2}
\end{eqnarray*}
for all $0\leq k\leq m-2$. Thanks to Lemma 4.3.2 of \cite{20}, we have
\begin{eqnarray*}
|\mathcal{L}^k_{\mathcal{A}}\psi_{\xi}|\lesssim r^{-1}\tau_{-}^{-1/2}\cdot\left\{\int_{\hat{\mathcal{C}}^t_v}\sum\limits_{i=0}^2|\mathcal{L}^{k+i}_{\mathcal{A}}\psi_{\xi}|^2\right\}^{1/2}
\end{eqnarray*}
and
\begin{eqnarray*}
|\mathcal{L}^k_{\mathcal{A}}\chi_{\iota}|\lesssim r^{-1}\tau_{-}^{-1/2}\cdot\left\{\int_{\hat{\mathcal{C}}^t_v}\sum\limits_{i=0}^2|\mathcal{L}^{k+i}_{\mathcal{A}}\chi_{\iota}|^2\right\}^{1/2}.
\end{eqnarray*}
As a consequence, for all $0\leq k\leq m-2$, the readers are able to gain
\begin{eqnarray*}
|(\mathcal{L}^k_{\mathcal{A}}\psi)_{\eta}|\lesssim\varepsilon_0\cdot r^{-3/2}\s\mbox{and}\s |(\mathcal{L}^k_{\mathcal{A}}\chi)_{\theta}|\lesssim\varepsilon_0\cdot r^{-3/2}
\end{eqnarray*}

\begin{eqnarray*}
|(\mathcal{L}^k_{\mathcal{A}}\psi)_{\xi}|\lesssim\varepsilon_0\cdot r^{-1}\tau_{-}^{-1/2}\s\mbox{and}\s |(\mathcal{L}^k_{\mathcal{A}}\chi)_{\iota}|\lesssim\varepsilon_0\cdot r^{-1}\tau_{-}^{-1/2}.
\end{eqnarray*}
Combining the above inequality leads to
\begin{eqnarray*}
|\mathcal{L}^k_{\mathcal{A}}\psi|\lesssim\varepsilon_0\cdot r^{-1}\max\{r^{-1/2},\tau_{-}^{-1/2}\}\s\mbox{and}\s |\mathcal{L}^k_{\mathcal{A}}\chi|\lesssim\varepsilon_0\cdot r^{-1}\max\{r^{-1/2},\tau_{-}^{-1/2}\},
\end{eqnarray*}
which is equivalent to
\begin{eqnarray*}
|\mathcal{L}^k_{\mathcal{A}}\Phi|\lesssim\varepsilon_0\cdot r^{-1}\max\{r^{-1/2},\tau_{-}^{-1/2}\}
\end{eqnarray*}
implying
\begin{eqnarray*}
|D^k\Phi|\lesssim\varepsilon_0\cdot r^{-1}\tau_{-}^{-k}\max\{r^{-1/2},\tau_{-}^{-1/2}\}.
\end{eqnarray*}
$\hfill\qedsymbol$

\subsection{Bound of the energy at the beginning with respect to initial data}\label{subsection3.4}
Given $0\leq k\leq m$, due to the structure of $Z_p$($Z_p\in\mathcal{A}$ and $0\leq p\leq k$) we have
\begin{eqnarray*}
(Z_0\cdots Z_k\Phi)(0,x)=\sum\limits_{|\alpha|+i+j\leq k}a^i(x)\cdot(\partial^{\alpha}_x\partial^j_t\Phi)(0,x)
\end{eqnarray*}
for some polynomials $a^i$ of degree $i$.

The nonlinear Dirac system we consider can be rewritten as
\begin{eqnarray}\label{9}
\partial_t\Phi+\gamma^0\gamma^j\partial_j\Phi=\sqrt{-1}\mathcal{N}(\Phi,\Phi)\gamma^0\Phi.
\end{eqnarray}
In the next, employing (\ref{9}) to transform time derivatives into spatial derivatives we write $(Z_0\cdots Z_k\Phi)(0,x)$ as
$$
\sum\limits_{i+|\alpha|\leq k}b^i(x)\cdot(\partial_x^{\alpha}\Phi)(0,x),
$$
where $b^i(x)$ is no more a polynomial, but a combination of $a^i(x)$ with $\partial_x^{\alpha}\partial_t^jA(0,x)$, $\partial_x^{\alpha}\partial_t^jB(0,x)$, $\partial_x^{\alpha}\partial_t^jC(0,x)$ and $\partial_x^{\alpha}\partial_t^jD(0,x)$. Because $\mathcal{N}$ is $C^m$-rapidly decreasing, we have
$$
(Z_0\cdots Z_k\Phi)(0,x)\lesssim\sum\limits_{i+|\alpha|\leq k}(1+|x|)^i|\partial^{\alpha}_x\Phi_0|.
$$
Hence $E(0)\lesssim||\Phi_0||^2_{H^m}$ and $E(0)+F(0)\lesssim||\Phi_0||^2_{H^m}$, since $F(0)=0$.

This completes the proof.
$\hfill\qedsymbol$

\end{document}